\documentclass[10pt]{amsart}
\usepackage{amssymb, amscd, amsmath, amsthm, epsf, epsfig, latexsym}
\usepackage{pb-diagram}
\def\zz{{\bf Z}}
 
\def\qq{{\bf Q}}

\def\cs{\mathop{\#}}

\def\calg{\mathcal{G}}
\def\calc{\mathcal{C}}

\newcommand{\fig}[2] { \includegraphics[scale=#1]{#2} }

\newtheorem{theorem}{Theorem}

\newtheorem{corollary}[theorem]{Corollary}

\theoremstyle{definition}
\newtheorem{definition}[theorem]{Definition}

\numberwithin{equation}{section}

\begin{document}

\title[Invariants of knots with slice Bing doubles]
{Algebraic and Heegaard-Floer invariants of knots with slice Bing doubles}

\author{Jae Choon Cha}
\address{Information and Communications University,
 Munji-dong, Yuseong-gu, Daejeon 305-732, Korea}
\email{jccha@icu.ac.kr}
\author{Charles Livingston}
\address{Department of Mathematics\\ Indiana University\\ Bloomington, Indiana 47405,  U.S.A.} 
\email{livingst@indiana.edu}
\author{Daniel Ruberman}
\address{Department of Mathematics, M050\\ Brandeis University\\ Waltham, Massachusetts 02454, U.S.A.}
\email{ruberman@brandeis.edu}
\thanks{The first author was partially supported by a grant from the KRF.  The second and third authors were partially supported by NSF Grants 0406934 and 0505605, respectively.}

%%%%%%%ABSTRACT%%%%%%%%%%%%%%

\begin{abstract}  If the  Bing double of a knot $K$ is   slice, then $K$ is algebraically slice.  In addition the Heegaard--Floer concordance invariants $\tau$, developed by Ozsv\'ath-Szab\'o, and $\delta$, developed by Manolescu and Owens, vanish on $K$.
 \end{abstract}

\maketitle

%%%%%%%SECTION%%%%%%%%%%%%%%

\begin{figure}[h]
\begin{center}
 \fig{.4}{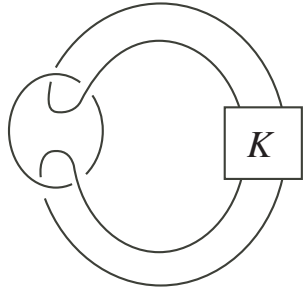} 
\end{center}
\caption{Bing doubling,  $B(K)$} \label{bingdoublefig}
\end{figure}

For a knot $K \subset S^3$, the Bing double, denoted $B(K)$, is the two component link  illustrated schematically in Figure~\ref{bingdoublefig}.  Within the box the two strands run parallel along a diagram for $K$, and for this to be well-defined, independent of the choice of diagram of $K$, the   strands are twisted so that their algebraic crossing number within the box is zero.  

A link in $S^3$ is called {\it slice} if its components bound disjoint locally flat  disks in the 4--ball, a notion that is sometimes called {\it strongly slice}. The slicing of links formed via the construction of Bing doubling has been a focus of  recent research, in part because of its connection to topological surgery in dimension four; see for instance~\cite{Kru1, Kru2, Kru3}.   Specifically, Bing doubling plays an essential role in the study of link theory, and especially in link concordance~\cite{Hab, KruTei, SchTei}.

 Harvey~\cite{Har} and   Teichner (unpublished) proved that if    $B(K)$ is slice, then the integral of the signature function of $K$ over the unit circle is zero.  Cimasoni~\cite{Cim} recently extended this in the case that $B(K)$ is boundary slice (that is, the slice disks along with Seifert surfaces for the components of $B(K)$ bound disjoint embedded 3--manifolds in the 4--ball), showing that this added assumption implies then $K$ is algebraically slice.    Cha~\cite{Cha1} strengthened the Harvey and Teichner result, showing that if $B(K)$ is slice, then the signature function is identically zero.  He also showed that some knots with vanishing signature function, such as the figure eight knot, have nonslice Bing doubles.  Cochran, Harvey and Leidy~\cite{HarOWR} have announced a proof that for certain algebraically slice knots $K$, $B(K)$ is not slice.  We prove the following:

\begin{theorem}\label{mainthm}   If $B(K)$ is slice, then $K$ is algebraically slice. \end{theorem}

As a corollary, this implies that the Arf invariant of $K$ is trivial.  This question alone had been a subject of  research interest.

Cimasoni also applied the Rasmussen link invariant~\cite{BelWeh, Ras} to find obstructions to a Bing double being smoothly slice  and used this to  show that if $K$ has Thurston-Bennequin invariant TB($K) \ge 0$ then $B(K)$ is not smoothly slice.  We show that if the Ozsv\'ath-Szab\'o invariant $\tau(K)$ (see~\cite{OzsSza}) or the Manolescu-Owens invariant $\delta(K)$ (see~\cite{ManOwe})   is nonzero, then $B(K)$ is not smoothy slice;  Cimasoni's result concerning TB($K)$ follows from this as well, using the connection between $\tau$ and TB, proved in~\cite{Liv}.

We have written this paper to make it as self-contained as possible.  But it is valuable to view it from the more general perspective of rational knot concordance.  This places the work in its historical context, provides a more general perspective, and provides more concise arguments built upon deeper theory.  In a final section we summarize this approach.

We thank David Cimasoni for discussing the topic of Bing doubling with us, and for his careful reading of an initial draft of this paper.

\section{Review of concordance and algebraic concordance}

Let $R \subset \qq$ be a subring of the rationals.  In this paper $R$ will be either $\zz$, $\qq$, or $\zz_{(2)}$, $\zz$ localized at 2.  This last ring is simply the set of rationals that can be written as   $a/b$ with $b$ odd.   Recall that for a space $X$, $H_i(X,\zz_{(2)}) = 0$ if and only if $H_i(X,\zz/2\zz) = 0$. 

The  Blanchfield pairing of a knot  $K \subset \Sigma^3$, where $\Sigma^3$ is an $R$-homology 3--sphere ($H_*(\Sigma^3, R) \cong H_*(S^3, R) $), arises as follows.  Up to sign   there is a unique surjection of $H_1(\Sigma^3 - K, \zz$) to $\zz$,  and this induces an infinite cyclic cover,  say~$X$.  We have that $H_1(X, R)$  is a torsion $R[t, t^{-1}]$--module and there is a nonsingular Hermitian  $S^{-1}_R R[t,t^{-1}]/R[t, t^{-1}]$--valued linking pairing on this module which is called the ``Blanchfield pairing.'' The Witt class of this pairing is 
denoted $W_R(K)$.  Here the Witt group is $L(R[t,t^{-1}],S_R)$, where $S_R = \{p(t) \in R[t,t^{-1}] \ |\ p(1) \text{\ is a unit in\ } R\}$.  This is the Witt group  of pairs $(H,\beta)$ where $H$ is an $R[t,t^{-1}]$--module which is $S_R$--torsion and $\beta$ is a 
nonsingular Hermitian   form  with values in $S^{-1}_R R[t, t^{-1}]/R[t, t^{-1}]$.   It is a theorem~\cite{Lev2} that the natural maps $L(\zz[t,t^{-1}],S_\zz) \to  L(R[t,t^{-1}],S_R)  \to  L(\qq[t,t^{-1}],S_\qq)$ are injections.   

\begin{definition} We say that $K$ is {\it algebraically slice}  if $W_R(K)$ is Witt trivial.  \end{definition}

\begin{theorem}\label{R2Zthm} If $K \subset S^3$ bounds a slice disk in a $R$--homology ball $Y^4$ with $R = \zz \text{\ or\ } \zz_{(2)}$, then $W_\zz(K) = 0$. 
\end{theorem}

A general discussion of this result is contained in~\cite[Section 2.2, 4.4]{Cha2}.  We give a brief proof in Section~\ref{proofsection}.  This statement is not true if $\zz_{(2)}$ is replaced with $\zz_{(p)}$ for $p$ odd; Cochran observed that the figure eight knot is slice in a rational homology ball (see~\cite{Cha2} for an extended discussion of such examples).

\section{Companionship and Blanchfield pairings}

Let $K_1 \cup U \subset S^3$ be a link with $U$ unknotted and linking number $w$,    let $K_2 \subset S^3$ be a knot, and let $N(K)$ denote the interior of a tubular neighborhood of $K$.  The union $(S^3 - N(U)) \cup (S^3 - N(K_2))$, where the union identifies the peripheral tori via a map that interchanges longitudes and meridians, is homeomorphic to $S^3$ and the image of $K_1$ in this union will be denoted ${K_2}( {K_1}, U)$.  In effect, $K_2(K_1,U)$ is the knot that results by tying the strands of $K_1$ that run through a disk bounded by $U$ into the knot $K_2$. We will see examples of this in the next section.  

It is proved in~\cite{LivMel} that the Blanchfield pairing of $K_2(K_1,U)$ is determined by those of $K_1$ and  $K_2$ as follows.  The homology of the infinite cyclic cover of  $K_i$ (that is, the {\it Alexander module}) is presented by a matrix $A_i(t)$ and the Blanchfield pairing is given by a matrix $B_i(t)$.  The main result of~\cite{LivMel} (see also~\cite{Kea}) states:

\begin{theorem}\label{companthm} The Alexander module and  Blanchfield pairing of $K_2(K_1,U)$ are presented by $A_1(t) \oplus A_2(t^w)$ and $B_1(t) \oplus B_2(t^w)$.
\end{theorem}

As a special case, if $K_1 \cup U$ is the $(2,2n)$--torus link (so that $K_1$ is also an unknot), $K_2(K_1,U)$ is simply the $(n,1)$--cable of $K_2$.  In general, letting $J^{(n)}$ denote the $(n,1)$--cable of a knot $J$ we have:

\begin{corollary} If the Alexander module and  Blanchfield pairing of $J$ are presented by $A(t)$ and $B(t)$, then the Alexander module and  Blanchfield pairing of $J^{(n)}$ are presented by $A(t^n)$ and $B(t^n)$.
\end{corollary}

\subsection*{The induced map on $L(\zz[t,t^{-1}],S_\zz)$}

Any class $W \in L(\zz[t,t^{-1}],S_\zz) $ is represented by a knot, $W = W_\zz(K)$.  The map that sends $K$ to $K^{(n)}$ induces a map $\phi_n :L(\zz[t,t^{-1}],S_\zz) \to  L(\zz[t,t^{-1}],S_\zz)$. (This map can also be defined via the map that sends $(A(t), B(t))$ to $(A(t^n), B(t^n))$.)

\begin{theorem}\label{alloddnthm}
  \mbox{}
  \begin{enumerate}
  \item $\phi_n\colon L(\zz[t,t^{-1}],S_\zz) \to L(\zz[t,t^{-1}],S_\zz)$ as
    described above is well-defined.
  \item For $n$ odd, $\phi_n$ is injective.
  \item If $\phi_n(W) = W$ for all odd $n$, then $W = 0$.
  \end{enumerate}
\end{theorem}

Theorem~\ref{alloddnthm} (2) seems to be folklore; the earliest
published reference we find is~\cite{CocOrr}, where it is stated
without indication of proof.  See also \cite[Section 4,4]{Cha2}.
(3)~is essentially a result of Kawauchi~\cite{Kaw}; we give a brief
proof in Section~\ref{proofsection}.

\section{Proof of Theorem~\ref{mainthm}: Algebraic Sliceness}\label{sectionpf}

\begin{figure}[h]
\begin{center}
\fig{.4}{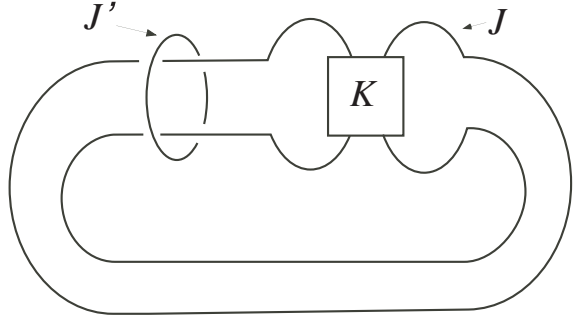}  
\end{center}
\caption{Alternative diagram,  $B(K)$} \label{bingdouble2}
\end{figure}

\begin{figure}[h]
\begin{center}
  \fig{.6}{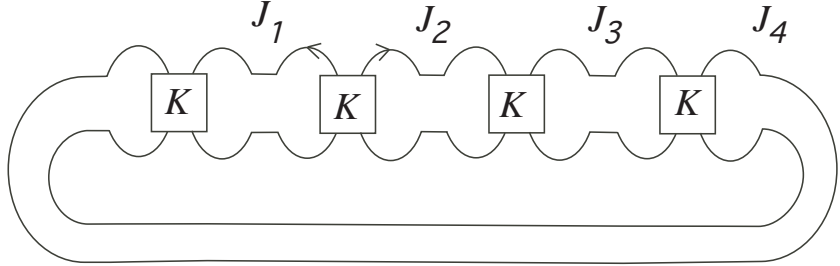} 
\end{center}
\caption{Four-fold covering link} \label{figurecover}
\end{figure}

We suppose here that  $B(K)$ is slice.  Figure~\ref{bingdouble2} gives an alternative diagram of $B(K)$, labeling its two components $J$ and $J'$.  The 4--fold cyclic branched cover of $B^4$ branched over the slice disk for $J'$ is a $\zz_{(2)}$--homology ball with boundary $S^3$ (since $J'$ is unknotted), which we denote $Y^4$.  (That is, $\tilde{H}_*(Y^4, \zz/ 2\zz) = 0$.)  The preimage of $J$ is a link of four components: $\tilde{J} = \{ J_1, J_2, J_3, J_4\}$ as illustrated in the Figure~\ref{figurecover}.  Notice that the slice disk for $J$ lifts to give a   slicing of $\tilde{J}$ in $Y^4$.

We consider only the link $\{J_1 , J_2\}$ which is   slice in $Y^4$, and which we now orient as shown in Figure~\ref{figurecover}.  It follows that for any $p$ and $q$, the   cable knots $(J_1)^{(p)}$ and $(J_2)^{(q)}$ are slice in $Y^4$.  (This depends on having initially taken the 0--framed Bing double, so that parallel push-offs of the $J_i$ are preimages of parallel push-offs of $J$, bounding disjoint push-offs of the slice disk.)  We now denote by   $J(p,q)$ the band connected sum $(J_1)^{(p)} \cs_b (J_2)^{(q)}$. This knot depends on the choice of band $b$; we can select $b$ so that it misses the box labeled ``$K$'' in the diagram and so that in the case that $K$ is unknot, the resulting knot $(J_1)^{(p)} \cs_b (J_2)^{(q)}$ is trivial. Clearly, as the band connected sum of slice knots, $J(p,q)$ is slice.  (This  construction is a special case of building links from covering links, as described in~\cite{ChaKo,CocOrr}.)

\begin{theorem}\label{pqthm} For the knot $J(p,q)$ constructed above, 
$$W_\zz(J(p,q)) =    \phi_p W_\zz (K) + \phi_{p+q} W_\zz (K) + \phi_q W_\zz (K) .$$
\end{theorem}

\begin{proof}
Let $J_0(p,q)$ denote the knot built as above, only with $K$ the unknot.  As noted above, the band $b$ is chosen so that $J_0(p,q)$ is an unknot.  Then $J(p,q)$ is built from $J_0(p,q)$ by three successive companionship constructions, with winding numbers $p$, $p+q$, and $q$.  Applying Theorem~\ref{companthm} to compute the Blanchfield pairing gives the desired result.  
\end{proof}

\noindent{\bf Proof of Theorem~\ref{mainthm}}

\begin{proof} Since $J(p,q)$ is slice in $Y^4$, by Theorem~\ref{R2Zthm}  we have $W_\zz(J(p,q)) = 0$, and by Theorem~\ref{pqthm}
$$\phi_p W_\zz (K) + \phi_{p+q} W_\zz (K) +\phi_q  W_\zz (K) = 0 $$ for all $p$ and $q$.

Letting $p=1$, this gives $$\phi_1 W_\zz (K) + \phi_{q+1} W_\zz (K) +\phi_q  W_\zz (K) = 0$$ for all $q$.  Replacing $q$ with $q-1$ gives $$\phi_1 W_\zz (K) + \phi_{q} W_\zz (K) +\phi_{q-1}  W_\zz (K) = 0$$ for all $q$.
Combining these gives $$  \phi_{q-1} W_\zz (K)  = \phi_{q+1}  W_\zz (K) $$ for all $q$.  In particular, $$  \phi_{q} W_\zz (K)  = \phi_{1}  W_\zz (K) = W_\zz(K) $$ for all odd $q$.  Finally, by Theorem~\ref{alloddnthm}, $W_\zz(K) = 0$.
\end{proof}

\section{Heegaard Floer Invariants}
As observed in the previous section,  if the Bing double of $K$ is slice, then $K \cs K^r$ is smoothly slice in a rational ball.  It follows from~\cite{OzsSza} that for  the  Ozsv\'ath-Szab\'o invariant $\tau$, if $J$ is smoothly  slice in a rational homology ball, then $\tau(J) = 0$.  Thus, in the present situation,  $\tau(K \cs K^r) = 0$.  But $\tau$ is an  additive $\zz$--valued invariant and does not detect knot orientation, so $\tau(K) = 0$.

The Manolescu-Owens invariant~\cite{ManOwe}, $\delta$, of a knot $J$ is defined to be the correction term of the Heegaard-Floer homology of the 2--fold branched cover of $S^3$ branched over $J$.   In our setting, $K \cs K^r$ is smoothly slice in a $\zz_{(2)}$--homology ball.  Thus, the 2--fold branched cover of $S^3$ branched over $K \cs K^r$   bounds a smooth $\zz_{(2)}$--homology ball.  It follows that the correction term for this cover is 0.  Again using the additivity of $\delta$, it follows that $\delta(K) = 0$.

\section{Proofs of theorems~\ref{R2Zthm} and~\ref{alloddnthm} }\label{proofsection}

\noindent{\bf Theorem~\ref{R2Zthm}.} {\it  If $K \subset S^3$ bounds a slice disk in an $R$--homology ball $Y^4$ with $R = \zz \text{\ or\ } \zz_{(2)}$, then $W_\zz(K) = 0$. }

 \begin{proof}  Let $K$ bound a slice disk $D \subset Y^4$.  Then $D$ represents an element in  $H_2(Y^4, \partial Y^4)$.  For every positive $n$,  the $(n,1)$--cable of $K$, $K^{(n)}$, bounds a disk $D^{(n)}$ built from $ n$ parallel copies of $D$.  Since $H_2(Y^4, \partial Y^4)$ has odd order, for some odd $n$, $D^{(n)}$ represents $0 \in  H_2(Y^4, \partial Y^4)$.  It follows that $K^{(n)}$ is boundary slice in $Y^4$, and thus the standard  argument  (that for knots in $S^3 = \partial B^4$, slice implies algebraically slice) can be applied and we see that $K^{(n)}$ is algebraically slice.  That is, $\phi_n W_\zz(K) = 0$.  By Theorem~\ref{alloddnthm}, to be proved next, this implies $W_\zz(K) = 0$, as desired.
 \end{proof}

\noindent{\bf Theorem~\ref{alloddnthm}.} {\it
  \begin{enumerate}
  \item The map $\phi_n :L(\zz[t,t^{-1}],S_\zz) \to
    L(\zz[t,t^{-1}],S_\zz)$ induced by $K \to K^{(n)}$ is a
    well-defined homomorphism.
  \item For $n$ odd, $\phi_n$ is injective.
  \item If $\phi_n(W) = W$ for all odd $n$, then $W = 0$.
  \end{enumerate}
}

\begin{proof}[Proof of (1)]
The map $\phi_n$ can be described as follows.  Any  given class $W \in  L(\zz[t,t^{-1}],S_\zz)$ can be represented as $W = W_\zz(K)$ for some knot $K$.  Then $\phi_nW = W_\zz ({K}^{(n)})$. To see the independence on the choice of $K$, suppose that  $W_\zz(K) = W_\zz(K_1)$ for some knot $K_1$.  Then $K \# -K_1$ is algebraically slice, so $W_\zz(K \# -K_1) = 0$.    There is a slice knot $L$ with the same Blanchfield pairing as $K \# - K_1$.  Since $L$ is slice,   $L^{(n)}$ is slice, and so has trivial Blanchfield pairing.  But according to~\cite{LivMel}, this pairing is determined by that of $L$, so the Blanchfield pairing of $(K \# -K_1)^{(n)}$ is trivial.  Again applying the formula of~\cite{LivMel}, this pairing is the same as the direct sum of the pairings for $K^{(n)}$ and $-{K_0}^{(n)}$.  Thus, $W_\zz(K^{(n)}) = W_\zz({K_0}^{(n)})$, as desired, and the map $\phi_n$ is well-defined.    Also, as just mentioned, according to the formula of~\cite{LivMel}, $\phi_n$ is additive.
\end{proof}

\begin{proof}[Proof of (2)]
  To begin the argument, note first that via composition we only need
  to prove the injectivity for $n$ an odd prime, $p$.

Suppose that a knot $K \subset S^3$ is such that $\phi_p W_\zz (K)  =0$;  that is, $K^{(p)}$ is algebraically slice.  Let $\tilde{K}^{(p)}$ denote the preimage of $K^{(p)}$ 
in   $\Sigma^3$,  the $p$--fold branched cover of $S^3$ over $K^{(p)}$.  We observe that $\Sigma^3$ is a $\zz$--homology sphere. Perhaps the easiest way of seeing this is via 
Fox's formula for the order of the homology of a $p$--fold branched cover of a knot $K$:   $|H_1(\Sigma^3)| = | \prod_{i=1}^{i=p-1} \Delta_K(\zeta^i) |$, where $\zeta$ is a $p$--
root of unity.  Since $\Delta_{K^{(p)}}(t) = \Delta_K (t^p)$ and $\Delta_K(1) = \pm 1$, we have in the present case that $ | H_1(\Sigma^3)| = 1$.

According to a formula given in Theorem~2.1 of~\cite{ChaKo},
a Seifert matrix $\tilde{A}$ of $\tilde{K}^{(p)}$ is determined by a
Seifert matrix of $K_p$; for a Seifert matrix $A$ of $K^{(p)}$, let
$\Gamma=(A-A^T)^{-1}A$.  Then $\Gamma^p - (\Gamma-I)^p$ is
nonsingular, and the matrix $\tilde A$ given by
\[
\tilde A = A-A^T\big(\Gamma^{p-1}-(\Gamma-I)^{p-1}\big)\big(\Gamma^p -
(\Gamma-I)^p)^{-1}\Gamma
\]
is a Seifert matrix for $\tilde{K}^{(p)}$.  (This matrix may be rational; according to~\cite{Cha2} it determines the class $W_\qq(K)$.) Furthermore, it is easily seen that the algebraic cobordism class of $A$ determines that
of~$\tilde A$.  Thus, the Blanchfield pairing for $\tilde{K}^{(p)}$ is
Witt trivial in $L(\qq[t,t^{-1}],S_\qq)$, and hence is trivial in $L(\zz[t,t^{-1}],S_\zz)$ also.

The Blanchfield pairing for $\tilde{K}^{(p)}$ can also be computed
geometrically as follows.  Let $E$ be the exterior of $K$ and $X$ be
its infinite cyclic cover.  The infinite cyclic cover $\tilde X$ of
the exterior of $\tilde{K}^{(p)}$ is built from $p$ disjoint copies of
$X$, with the deck transformation acting individually on each copy of
$X$ via the original $K$ deck transformation.  More precisely, $\tilde
X$ is described as follows.  Regarding $K^{(p)}$ as a satellite knot, its
exterior consists of $E$ and the exterior, say $E'$, of the Hopf link
with one component replaced by its $(p,1)$-cable.  Considering the
preimage of this decomposition, one can see that $\tilde X$ is the
union of $pX$ and an infinite cyclic cover $X'$ of $E'$.  We have that $X'$ is the
exterior of $p$ disjoint long arcs in $D^2\times \mathbf{R}$, and so
$H_1(X')$ is a free abelian group generated by meridians of the long
arcs.  Therefore, from a Mayer-Vieotoris argument, it follows that
$H_1(\tilde X)=H_1(X)^p$.  Hence, $0=W_\zz( \tilde{K}^{(p)} )=pW_\zz(K)$.

Since the Witt group contains no odd torsion~\cite{Lev2}, it follows that $W_\qq(K)$
must have been Witt trivial.  This shows that $K$ is algebraically
slice.
\end{proof}

\begin{proof}[Proof of (3)]
 We wish to show that if $W \in L(\zz[t,t^{-1}],S_\zz)$ and $\phi_n(W )
= W $ for all odd $n >1$, then $W = 0$.  This is essentially a result
of Kawauchi, proved in~\cite{Kaw}.  We give a short self-contained
proof here.

Without loss of generality we assume that $W = W_\zz(K)$ for some knot
$K \subset S^3$.  We now switch to rational coefficients and consider
the Witt class $W_\qq(K)$ represented by the Blanchfield pairing
$\beta$ defined on $H=H_1(X;\qq)$ where $X$ is the infinite cyclic
cover of $S^3-K$.  The order of $H$ as a $\qq[t, t^{-1}]$--module is
$\Delta_K(t)$, the Alexander polynomial of $K$.

Suppose that $W_\qq(K) \ne 0$.  Then for some symmetric irreducible
polynomial $\lambda(t)$, the restriction of $\beta$ to the
$\lambda(t)$--primary part of $H$ is nontrivial in the Witt group.
Such $\lambda(t)$ always divides $\Delta_K(t)$.  Since $W_\qq(K) =
W_\qq(K_n)$ for $n$ odd, we also have that $\lambda(t)$ divides
$\Delta_{K_n}(t)$.  Since $\Delta_{K_n}(t)$ is equal to
$\Delta_K(t^n)$ up to units in $\qq[t,t^{-1}]$ (originally proved by Seifert~\cite{Sei}), $\lambda(t)$ divides $\Delta_K(t^n)$ for all odd~$n$.

Let $\alpha$ be a root of $\lambda(t)$.  Then for all odd $n$,
$\alpha^n$ is a root of $\Delta_K(t)$. But $\Delta_K(t)$ has only a
finite number of roots, so $\alpha^n = \alpha^m$ for some $n \ne m$.
It follows that $\alpha$ is a root of unity.  Since $t-1$ does not
divide any Alexander polynomial, $\alpha\ne 1$.  So $\lambda(t) =
\Phi_d(t)$ for some $d\ge 2$, where $\Phi_d(t)$ denotes the
$d$--cyclotomic polynomial.  For some prime power $p^k>1$ and odd $r$,
$\alpha^r$ is a primitive $p^k$-th root of unity; for, writing
$d=2^ap_1^{a_1}\cdots p_r^{a_r}$ where the $p_i$ are distinct odd
primes, if $a\ne 0$ then $\alpha^{d/2^a}$ is a primitive $2^a$-th root
of unity; if $a = 0$ then we may assume $a_1\ne 0$, and
$\alpha^{d/p_1^{a_1}}$ is a primitive $p_1^{a_1}$-th root of unity.
We now have that $\Delta_K(t)$ has as a root a prime power root of
unity.  This implies that $\Phi_{p^k}(t)$ divides $\Delta_K(t)$.  By
Gauss' lemma, $\Delta_K(t)=\Phi_{p^k}(t) f(t)$ for some $f(t)$ with
integral coefficients. But $\Phi_{p^k}(1) = p$ and $\Delta_K(1) = \pm
1$, so we have a contradiction.
\end{proof}

\section{Rational knot concordance}

In this section we want to discuss the previous work from the
viewpoint of general rational knot concordance.  Actually this
approach led us to the proof of Theorem~\ref{mainthm} presented above.
Our main reference is the monograph by the first author,~\cite{Cha2}.

For a subring $R \subset \qq$ the concordance group, $\calc^R$, is
built by considering knots in $R$--homology 3--spheres, with
concordances taking place in $R$--homology cobordisms. 

Our link $J(p,q)$ is a special case of the general construction of
forming a knot as the band sum of a {\it covering link}, that is, a
link formed as the union of components of the preimage of a link $L
\subset S^3$ in a $d$--fold branched cyclic covering of $S^3$ branched
over a component of~$L$ which is known to be a $\zz_{(p)}$-homology
sphere.  If $L$ is slice and $d=p^a$ is a prime power, then any knot
constructed in this way $L$ will be trivial in $\calc^{\zz_{(p)}}$.

Analogous to Levine's homomorphism of the (integral) knot concordance
group $\calc$ into the ``algebraic concordance group'' $\calg$, as
described in~\cite{Cha2} there is a natural group $\calg^{R}$ and a
homomorphism $\phi_R : \calc^R \to \calg^R$.  There is the following
commutative diagram.

\[
\begin{diagram}
  \node{\calc} \arrow{s}\arrow{e}
  \node{\calg} \arrow{s}
  \\
  \node{\calc^R} \arrow{e}
  \node{\calg^R}
\end{diagram}
\]

Consider a slice link $L$ with unknotted components.  Then its
covering links are always in~$S^3$.  In~\cite[Section 2.2, 4.4]{Cha2}
it is proved that for $R = \zz_{(2)}$ the map $\calg \to \calg^R$ is
injective.  It follows that any knot in $S^3$ that is built by banding
together components of a covering link of a slice link represents the
trivial element in $\calg$.  That is, it is algebraically slice.

Applied to $J(p,q)$, this implies (as was shown in
Section~\ref{sectionpf}) that $\phi_n W(K) = W(K)$ for all odd $n$.
Thus, Kawauchi's result (Theorem~\ref{alloddnthm}~(3)) can be applied
to conclude that $W(K) = 0$.

%%%BIBLIOGRAPHY%%%%%%%%

\newcommand{\etalchar}[1]{$^{#1}$} 

\end{document}